%% file: main.tex
\documentclass[11pt]{article}
\usepackage[margin=1.2in]{geometry}
\usepackage{amsfonts,amsmath,amssymb,amsthm,bm}
\usepackage{pifont}
\usepackage{pgfplots}
\usepackage{subcaption}
\usepackage{mathtools}
\usepackage{graphicx}
\usepackage{authblk}
\usepackage{algorithm}
\usepackage{algpseudocode}
\usepackage{subcaption}
\usepackage{caption}
\usepackage{tikz}
\usepackage{hyperref}

\title{Multi-stream Quickest Change Detection: Foundations and Recent Advances}

\author{Topi Halme}
\author{Visa Koivunen}
\affil{Aalto University, Finland}
\affil[ ]{\texttt{\{topi.halme, visa.koivunen\}@aalto.fi}}
\date{April 2026}

\newtheorem{Definition}{Definition}

\usepackage{subcaption}
\begin{document}
\maketitle
\begin{abstract}
This paper provides an overview of recent developments in quickest change detection (QCD) for high-dimensional multi-sensor systems, with an emphasis on settings involving structural constraints and limited sensing resources. Classical QCD methodologies, while well understood in low-dimensional and fully observed regimes, face significant challenges when extended to modern applications characterized by large-scale data, constrained sampling or communication, and heterogeneous signal structures.
We review key approaches for handling high dimensionality, including methods that exploit sparsity, and other forms of signal heterogeneity. Additionally, we discuss sampling constraints, where observations must be selected or acquired sequentially under resource limitations. Multi-stream applications can require making multiple detections, for example when detecting changes separately in different streams. The underlying assumptions on probability models, the types of changes taking place, commonly used decision-making criteria, performance indices, and error types are described. We also briefly discuss the application of machine learning in cases where the underlying probability models are not known or there is a need to select which sensors should monitor the phenomena because of the large scale of the system.
\end{abstract}


\include{shorts}

\section{Introduction}

In many modern applications, a fundamental problem is detecting and localizing rapid changes, anomalies or adversarial events in monitored processes or phenomena. Multiple data streams may be monitored simultaneously by using multiple devices or sensors in distinct locations.  Detecting a change in the statistical properties of a monitored process as quickly as possible is known as a \textit{quickest change detection} (QCD) problem \cite{POOR_BOOK,TARTAKOVSKY_BOOK}. The general objective in a QCD problem is to construct procedures that detect the change with as minimal delay as possible, while controlling the rate of false alarms. A false positive is considered to occur when one prematurely makes a decision that a change has taken place before it actually happened. Importantly, the inference must be done in real time, as in a variety of important applications the data is streaming, dynamic and generated sequentially.

Historically, the change detection problem was first studied in the context of industrial quality control during the 1920s and 1930s. The first algorithms are due to Shewhart \cite{SHEWHART_1925, SHEWHART_1930} and Page \cite{PAGE}. Shewhart proposed a scheme where an alarm is raised the first time an observation (or a function of it) falls outside prespecified limits. As such, the decision to stop or continue is only based on the most recent data point and does not utilize previous observations. The resulting detection rule is easy to implement, but the suboptimal use of all available data may make the detection of small changes very slow. Page \cite{PAGE} suggested using moving averages or cumulative sums of past statistics for more efficient detection. The algorithms of Shewhart and Page were not motivated by rigorous theoretical analysis, but by easy implementation and reasonable practical performance in quality control tasks. The first optimality considerations appeared in the 1950s and 60s after Wald's pioneering work on sequential analysis \cite{WALD_BOOK}. In a Bayesian context, the first optimal algorithms and theory are due to Shiryaev in a series of works in the 1960s \cite{SHIRYAEV_1963, SHIRYAEV_1965}. Optimality in a non-Bayesian, or minimax, sense was first studied by Lorden in 1971 \cite{LORDEN}. Since then, the general sequential change detection problem has been extended to many directions relevant to modern applications, including non-i.i.d. data models, uncertainty in the pre- and/or post-change distributions, high-dimensional problems, detection under sampling constraints, and decentralized systems, among others. Extensive general treatments of quickest change detection can be found in e.g. the books by Hadjiliadis and Poor \cite{POOR_BOOK} and Tartakovsky et al. \cite{TARTAKOVSKY_BOOK}.

In contemporary applications of QCD algorithms, a defining characteristic is a large number of data streams that must be monitored simultaneously in real time. A large number of devices equipped with sensors and communications capability may be observing streaming data continuously, for example smartphones, smart buildings, distributed radars and wireless base-stations, weather stations, seismic and environment monitoring, the Internet of Things, and sensors monitoring the state of our infrastructure (power grid, roads and bridges) to name a few. Depending on the application, the task may be to detect changes in individual streams separately. Alternatively, the goal may be to detect a single change that affects one or more of the streams simultaneously or with a small delay due to the physics of the phenomena, such as radio wave propagation. In both cases, the large number of streams poses additional, but different, difficulties and opportunities. In particular, this overview article highlights three central challenges that arise specifically in the multi-stream settings, and discusses solutions proposed in the literature. The challenges are:
\begin{enumerate}
    \item \textbf{Changes with high-dimensional unknown parameters:} In practical applications of sequential change detection, the distribution of the data after the change is usually not exactly known. Instead, the post-change distribution of the data may belong to a family of distributions that depends on an unknown parameter vector $\btheta$. For example, each component of $\btheta$ could represent the post-change mean parameter in different data streams, or $\btheta$ could denote the indices of sensors that have been affected by the change. If the number of streams is large, estimating $\btheta$ accurately becomes increasingly difficult, which in turn makes detecting a change more challengin.
    \item \textbf{Adaptive sampling under constraints:} When the number of sensors is large, it is often not feasible to receive data from all sensors at all times due to e.g. power and bandwidth constraints. Naturally, a tradeoff exists between the proportion of data observed and the ensuing detection delay. However, if the observer, such as a common fusion center, can adaptively choose which sensors it wants to receive observations from at a given time instance, the impact of reduced communication detection efficiency can be mitigated. Designing sampling policies that most effectively sample the sensors specifically for change detection purposes is therefore an important topic with lots of recent and ongoing research. 
    \item \textbf{Adjusting for multiple comparisons:} In many applications changes can happen at different times at different sensors, leading to many individual detections. In these cases, conventional QCD performance metrics need to be adjusted to take into account the simultaneous tests. Multiple Hypothesis Testing (MHT) provides a framework for controlling the error rates in that setting. 
    
\end{enumerate}

This paper presents a brief overview of underlying models and methods in multi-stream QCD. The underlying assumptions on probability models, the types of changes taking place, commonly used decision-making criteria, performance indices, and error types are described. We focus on large-scale settings where the number of streams and the dimension of the probability model parameter vectors may be large. We will also briefly discuss the application of machine learning in cases where the underlying probability models are not known or there is a need to select which sensors should monitor the phenomena because of the large scale of the system. We will describe some timely and emerging applications of multi-stream QCD. 

\noindent \textbf{Notation.}
Time is indexed by $n \in \mathbb{N}$. In the single-stream setting,
$X_1,X_2,\ldots$ denotes the observed data sequence, while in the multi-stream setting $X_n = [X_n^{(1)},\ldots,X_n^{(K)}]$ denotes the $K$-dimensional observation vector at time $n$,
where $K$ is the number of streams or sensors. The change-point is denoted by $\nu$, and in settings with streamwise changes $\nu^{(k)}$ denotes the change-point in stream $k$. Expectations and probabilities when change happens at time $\nu$
are denoted by $\mathbb{E}_\nu$ and $\mathbb{P}_\nu$, respectively, while $\mathbb{E}_\infty$ and $\mathbb{P}_\infty$ correspond to the no-change setting. When the post-change model depends on an unknown parameter $\btheta$, the corresponding expectation is denoted by $\mathbb{E}_\nu^\btheta$.
Pre- and post-change distributions are denoted by $f_0$ and $f_1$ in the simple setting, and more generally by $f(\cdot \mid \btheta)$ for a parametric family. 

\section{Quickest change detection: basic fundamentals}

In this section, we briefly introduce the standard QCD problem in the single-stream i.i.d. setting. The objective is to detect an abrupt change, if one occurs, in the underlying probability model of a data sequence $X_1,X_2,\ldots$ as quickly as possible. Initially, the data $X_1, \ldots, X_{\nu-1}$ is assumed to be generated according to a probability density $f_0$ which usually represents a normal operational state of the system. A common example is that the data observed by a sensor contains Gaussian noise only. At an unknown time $\nu$, the data distribution changes to some other distribution $f_1$, so that $X_{\nu}, X_{\nu+1},\ldots \sim f_1$. 

A change detection procedure is defined by a \emph{stopping time} random variable $T$ which indicates the time the sampling is stopped and a change is declared.
\begin{Definition}[Stopping time]\label{def:stoppingtime}
    Let $\calF_n = \sigma(X_1,\ldots,X_n)$ denote the sigma-algebra generated by the first $n$ observations.
    A random variable $T$ is a stopping time with respect to a sequence $X_1, X_2,\ldots$ if for each $n$, the event $\left\{T=n\right\} \in \calF_n$. Equivalently, the event $\left\{T=n\right\}$ is a function of only $(X_1,\ldots,X_n)$
\end{Definition}
A central modeling choice is whether the change-point $\nu$ is considered to be a random variable with a prior distribution, or an unknown deterministic quantity. The former yields the Bayesian approach, while the latter corresponds to a non-Bayesian, or minimax, framework. The two paradigms are introduced next. 

\subsection{Classical non-Bayesian objectives}
When $\nu$ is unknown, the usual QCD objective is to minimize some notion of average detection delay (ADD) while maintaining the average run length (ARL) to false alarm above a tolerable level $\gamma$. The ARL to false alarm is defined as 
\begin{equation}\label{eq: arl_definition}
    \ARL(T) = \Ex_\infty(T),
\end{equation}
meaning that when no change occurs, the constraint $\ARL(T) \geq \gamma$ implies that the test should take at least $\gamma$ samples on average before (false positive) detection. In other words, it is the expected number of observations to an alarm assuming that there is no change. Average detection delay is usually measured in the worst-case over values of $\nu$. Lorden \cite{LORDEN} proposed a ``worst-worst-case'' average detection delay (WADD) criterion, defined as
\begin{equation}\label{eq:Lorden_delay}
    \WADD(T) = \sup_{\nu \geq 1} \text{ess}~\!\sup \Ex_\nu[(T-\nu+1)^+ |\calF_{\nu-1}],
\end{equation}
where a supremum is first taken over all possible realizations of pre-change data $\calF_{\nu-1}$, and then all possible values of the change point. This leads to the following optimization problem
\begin{equation}\label{eq:Lorden_prob}
    \text{minimize}~~ \WADD(T) \quad \text{subject to } \ARL(T) \geq \gamma.
\end{equation}
In an alternate approach, Pollak \cite{POLLAK} suggested taking the average instead of the supremum over all possible pre-change observations. He proposed using conditional average detection delay (CADD), the expected number of samples needed to detect the change after it occurs as a criterion:
\begin{equation}\label{eq:CADD_def}
    \CADD(T) = \sup_{\nu \geq 1} \Ex_\nu(T-\nu | T\geq \nu).
\end{equation}
It is easy to see that $\WADD(T) \geq \CADD(T)$ for all $T$. 

Finding exactly optimal solutions to formulations like \eqref{eq:Lorden_prob} for all values of $\gamma$ is usually difficult except in few simple cases. For this reason, a common consideration in the literature is asymptotic optimality, usually studied in the  limit $\gamma \to \infty$. In the i.i.d. setting, a universal lower bound for CADD of any procedure as $\gamma \to \infty$ was derived by Lai \cite{LAI} as
\begin{equation}\label{eq:asym_LB}
    \inf_{T : \ARL(T) \geq \gamma} \CADD(T) \geq \frac{\log \gamma}{I}(1+o(1)),
\end{equation}
where $I$ is the Kullback-Leibler (KL) divergence (relative entropy) between distributions $f_1$ and $f_0$
\begin{equation}
    I = \int f_1(x) \log \frac{f_1(x)}{f_0(x)} dx.
\end{equation}
Earlier, a similar lower bound for WADD was derived by Lorden \cite{LORDEN}. A change detection procedure $T$ is said to be asymptotically (first-order) optimal if
\begin{equation}
    \frac{\CADD(T)}{\inf_{\tau : \ARL(\tau) \geq \gamma} \CADD(\tau)} \to 1 \quad \text{as } \gamma \to \infty.
\end{equation}

\subsection{Fundamental non-Bayesian QCD algorithms}

Next, we describe briefly two foundational algorithms for QCD in the non-Bayesian paradigm, i.e., the CuSum procedure and the Shiryaev-Roberts procedure.

\noindent\textbf{CuSum procedure. }The Cumulative Sum (CuSum) procedure was first introduced by Page \cite{PAGE}. The CuSum test tracks the drift of the cumulative log-likelihood-ratio sum $S_1,S_2,\ldots$:
\begin{equation}
    S_n = \sum_{j=1}^n \log \frac{f_1(X_j)}{f_0(X_j)}.
\end{equation}
It stops when the sum exceeds its past minimum by a sufficient amount. The test is defined by
\begin{align}\label{eq:cusum_def}
    W_n     & = S_n - \min_{k \leq n} S_k, \\
    \TCUSUM & = \inf\{n : W_n > b\}.
\end{align}
The CuSum statistic $W_n$ can alternatively be written as
\begin{equation}\label{eq: CuSum max form}
    W_n = \max_{k\leq n}\sum_{j = k}^n \log\frac{f_1(X_j)}{f_0(X_j)}
\end{equation}
and updated recursively using the formula
\begin{equation}\label{eq:cusum_recursion}
    W_n = \max\left\{0, W_{n-1} + \log\frac{f_1(X_n)}{f_0(X_n)}\right\}, \quad W_0 = 0.
\end{equation}
The recursive update makes the CuSum test efficient to implement in practice, and provides a low complexity method for detecting changes in streaming data. Lorden \cite{LORDEN} proved that as $\gamma \to \infty$, the CuSum procedure is asymptotically first-order optimal in terms of WADD. Later Moustakides \cite{MOUSTAKIDES} showed that CuSum with the stopping threshold $b$ chosen such that $\ARL(\TCUSUM) = \gamma$ is exactly optimal for all $\gamma$ in the sense of \eqref{eq:Lorden_prob}. Finding a $b$ such that the ARL constraint is satisfied with equality is generally not possible in closed form, but can be found e.g. numerically via Monte Carlo simulations. A simple lower bound for $\ARL(\TCUSUM)$ as a function of $b$ exists, and is given by (see e.g. \cite[Lemma 8.2.1.]{TARTAKOVSKY_BOOK})  
\begin{equation}\label{eq:arl_LB}
    \ARL(\TCUSUM) \geq e^b.
\end{equation}
Choosing $b = \log \gamma$ therefore guarantees $\ARL(\TCUSUM) \geq \gamma$.

\noindent\textbf{Shiryaev-Roberts procedure. }The Shiryaev-Roberts (SR) procedure is defined by the statistic
\begin{equation}\label{eq:SR_stat}
    R_n = \sum_{k=1}^n \prod_{j=k}^n \frac{f_1(X_j)}{f_0(X_j)}
\end{equation}
and the stopping time
\begin{equation}\label{eq:SR_def}
    \TSR = \inf\{n : R_n > b\}
\end{equation}
where $b$ is a threshold. The difference between the SR and CuSum statistics is that instead of taking the maximum over $k\leq n$ as in \eqref{eq: CuSum max form}, the SR statistic takes the sum of likelihood ratios over $k \leq n$. The SR statistic also has a recursive form, given by
\begin{equation}
    R_n = (1+R_{n-1})\frac{f_1(X_n)}{f_0(X_n)}, \quad R_0 = 0.
\end{equation}
The recursion stops and raises an alarm when $R_n$ exceeds the threshold $b$.
The procedure is also asymptotically optimal with respect to both WADD and CADD \cite{POLUNCHENKO_ON_OPTIMALITY}. The same ARL lower bound \eqref{eq:arl_LB} as for CuSum also applies to SR.

\subsection{Bayesian objective}\label{sec:bayes_object}
In some settings, it may be possible to consider the change-point $\nu$ as a random variable with a known distribution $\pi_n = \P(\nu = n)$. Optimal stopping times under this formulation were first studied by Shiryaev \cite{SHIRYAEV_1963} under the loss function 
\begin{equation}\label{eq:bayes_loss}
    B(T, \lambda) = \P^\pi(T < \nu) + \lambda\Ex^\pi[(T-\nu)^+],
\end{equation}
where the notation $\P^\pi, \Ex^\pi$ is used to highlight that $\nu \sim \pi$. With an appropriate choice of $\lambda$, the objective \eqref{eq:bayes_loss} is equivalent to the following constrained optimization problem:
\begin{equation}\label{eq:constrained_bayes}
    \text{minimize}~~ \Ex^\pi[(T-\nu)^+] \quad \text{subject to   }~ \P^\pi(T < \nu) \leq \alpha,
\end{equation}
i.e. minimizing the average detection delay under a constraint that the probability of false alarm (PFA) $\P^\pi(T < \nu)$ is smaller than some $\alpha$. When $\nu \sim \text{Geom}(\rho)$, Shiryaev \cite{SHIRYAEV_1963} showed that the exactly optimal test for \eqref{eq:bayes_loss} is a threshold test on the posterior probability of change
\begin{equation}\label{eq:posterior_prob}
    p_n = \P^\pi(\nu \leq n | \calF_n).
\end{equation}
The probability $p_n$ can be computed recursively as
\begin{equation}
    p_n = \frac{\Tilde{p}_{n-1}L(X_n)}{\Tilde{p}_{n-1}L(X_n) + (1-\Tilde{p}_{n-1})},
\end{equation}
where $\Tilde{p}_n = p_n + (1-p_n)\rho$ and $p_0 = \pi_0 = \rho$.
Shiryaev's test $T^S_b$ is then performed by comparing the probability $p_n$ to a stopping threshold $b > 0$, i.e.
\begin{equation}\label{eq:shiryaev_time}
    T^S_b = \inf\{n : p_n \geq b\}.
\end{equation}
As $\rho \to 0$, i.e. the prior distribution becomes increasingly uninformative, the limit of posterior odds $p_n/(1-p_n)$ converges to the Shiryaev-Roberts statistic \eqref{eq:SR_stat}.
As for CuSum and SR, finding a threshold $b$ such that the constraint in \eqref{eq:constrained_bayes} is satisfied with equality is usually only possible numerically. However, by the law of iterated expectation
\begin{equation}
    \P^\pi(T^S_b < \nu) = \Ex^\pi(\P^\pi(T^S_b < \nu | \calF_{T^S_b})) = \Ex^\pi(1-p_{T^S_b}) \leq 1-b,
\end{equation}
where inequality follows from the definition of $T^S_b$ \eqref{eq:shiryaev_time}. The probability of false alarm constraint in \eqref{eq:constrained_bayes} is then satisfied by setting $b = 1-\alpha$.

Shiryaev's test is also asymptotically optimal as false alarm probability $\alpha \to 0$ in a large class of non-i.i.d. models and for general prior distributions $\pi$ \cite{TARTAKOVSKY_GENERAL}. 

\subsection*{Comparison of Shiryaev's test and CuSum}

While the non-Bayesian and Bayesian procedures have very different notions of optimality, their empirical performance can be compared under common criteria. To illustrate this by a concrete example, we simulate a simple change-detection setting of a mean-shift from $0$ to $\mu = 0.5$ under Gaussian noise, so that $f_0 = \calN(0,1), f_1 = \calN(\mu, 1)$. Four algorithms are considered, the CuSum test \eqref{eq:cusum_def}($\TCUSUM$), SR procedure \eqref{eq:SR_def} ($\TSR$), and two versions of Shiryaev's test; one with parameter $\rho =1/1000$ ($T^{S}_{b,1000}$) and the other with $\rho = 1/50$ ($T^{S}_{b,50}$). For each procedure $T$, the detection threshold $b$ is set such that $\ARL(T) = \Ex_\infty(T) = 1000$ for purposes of comparison, where the values are found by 5000 Monte Carlo simulations. Note that the Shiryaev procedures are typically calibrated to have a prespecified PFA $\P^\pi(T<\nu)$ rather than ARL. PFA depends on $\pi$, the assumed prior distribution of the change-point $\nu$, which may not be available. The ARL guarantee holds regardless of $\nu$, even for the Bayesian procedures. 

\begin{figure}
\subfloat[\label{fig:simple_compare}]
{\includegraphics[width=7.0cm]{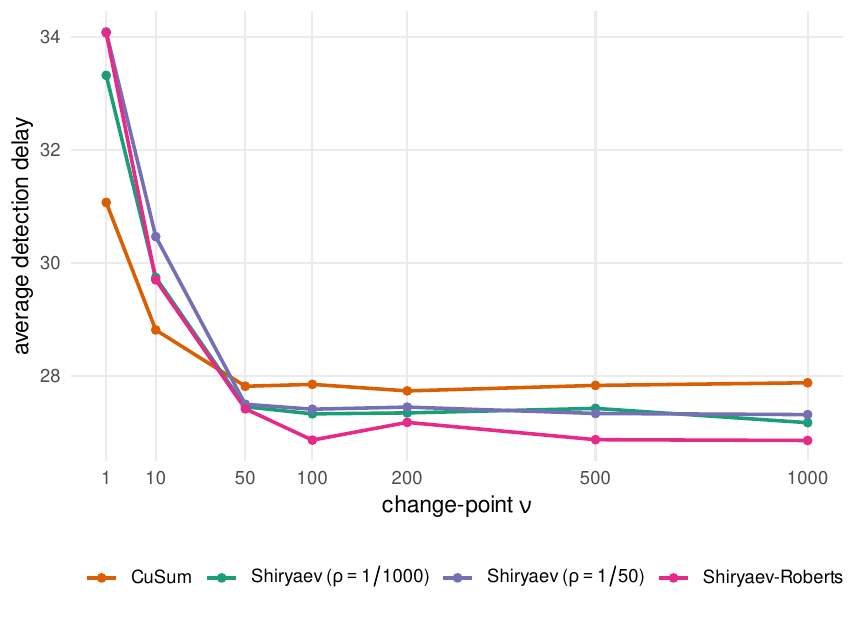}}
\hfill
\subfloat[\label{fig:arl_add}]
{\includegraphics[width=7.0cm]{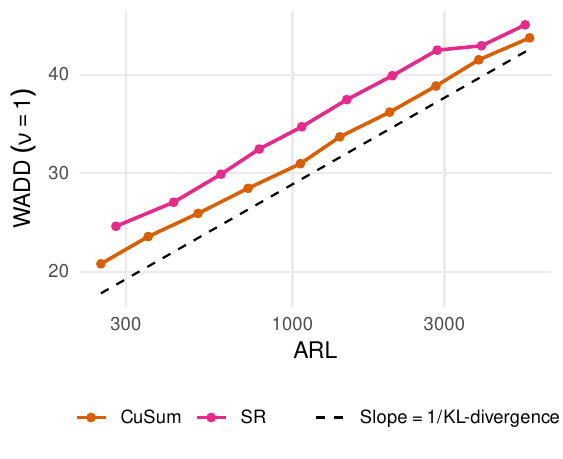}}

\caption{Change-point detection algorithms trade off between controlling false alarm rates (such as average run length (ARL) to false alarm) and average detection delay. \textbf{(a)} Average detection delay as a function of different values of the change-point $\nu$. CuSum procedure minimizes worst-case detection delay, but for many values of $\nu$, the delay of Shiryaev-Roberts or Bayesian Shiryaev tests (calibrated to have the same ARL) is smaller. \textbf{(b)} Worst-case average detection delay for CuSum and SR (which happens when $\nu =1$ as observed in the left panel) as a function of $\log($ARL$)$. For both algorithms, the detection delay grows linearly with respect to $\log$(ARL), and the slope is given by the reciprocal of KL-divergence, in accordance with asymptotic lower bound eq.~\eqref{eq:asym_LB}.}
\end{figure}

The average detection delay $\Ex_\nu(T- \nu | T\geq \nu)$ is evaluated for different values of $\nu$ in Figure~\ref{fig:simple_compare}) for all procedures. This delay metric is different than either $\CADD$ and $\WADD$, which include taking the maximum over $\nu$. For all tests, the largest detection delay occurs when the change happens early at time $\nu = 1$. This is intuitively explained by the fact that the test statistics that are typically updated using a recursive formula are guaranteed to be far away from the threshold if the change takes place in the very beginning of the monitoring process. The optimality of CuSum in terms of $\WADD$ is observed from the figure, since its maximum delay over $\nu$ is the smallest. On the other hand, for larger values of the change-point, the detection delay of Shiryaev-Roberts or the Bayesian procedures is smaller. The performance of the Shiryaev and SR procedures is similar, which is to be expected as in the limit $\rho \to 0$ the Shiryaev procedure converges to SR. In Figure~1b), the WADD of CuSum and SR is plotted as function of $\log(\ARL)$. The relationship is linear, and the slope given by the reciprocal of KL-divergence $(\mu^2/2)^{-1}$ matches the asymptotic lower bound~\eqref{eq:asym_LB}.

\subsection{Alternative performance criteria}

In single-stream QCD, minimizing some notion of average detection delay while controlling the ARL to false alarm (in non-Bayesian settings) or probability of false alarm (in Bayesian settings) is the most common performance objective. However, in some applications, other performance criteria may also be of interest. For instance, if $T$ represents the time when a machine part deemed to be faulty is replaced, it may make sense to penalize, for example, the average ``size of miss" $\Ex^\pi|T-\nu|$. More generally, the objective trades off between early and late detection 
\begin{equation}
    R(T, \lambda) = \Ex^\pi(\nu-T)^+ + \lambda\Ex^\pi(T-\nu)^+,
\end{equation} 
where $\lambda$ is a tunable parameter that controls the relative importance of the two criteria. Under such a penalty, replacing the part before it is actually faulty is penalized less if the replacement is done close to the eventual time of failure $\nu$. It turns out, that under some conditions, the Shiryaev test $T^S_b$ is also optimal for the risk $R(T, \lambda)$ \cite{Shiryaev2004}. 

In non-Bayesian settings, exactly optimal algorithms are known in relatively few cases besides the optimality of CuSum for Lorden's delay with i.i.d data. As an interesting example, Poor \cite{POOR_EXPONENTIAL} considered an exponential delay penalty of the form
\begin{equation}\label{eq:exponential}
    \dfrac{\beta^{(T-\nu+1)^+}-1}{\beta - 1}
\end{equation}
for some constant $\beta \not = 1$. If $\beta > 1$, each additional unit of delay is penalized in an exponentially increasing manner. It may be necessary to associate a higher cost on longer delays when the cost of undetected changes may be enormous. If $\beta < 1$, the cost of delay is saturating and sub-linear, with $\beta \to 1$ converging to the traditional linear detection penalty. For instance, in financial applications, a delayed detection of a shift in the market may incur an exponentially increasing penalty due to the compounding nature of financial gains and losses. Similarly, in applications where the change-point triggers a cascade of events (e.g. healthcare, power grid, communication networks, etc.), the effects of delay may be multiplicative, motivating an exponential penalty. Poor showed that a modified version of the CuSum procedure with recursion 
\begin{equation}\label{eq:exp_cusum_recursion}
    W_n = \max(0, W_{n-1}) + \ell(X_n) + \log \beta,  \quad W_0 = 0
\end{equation} 
minimizes the worst-case average of the exponential delay in \eqref{eq:exponential} under an ARL constraint. The modified CuSum statistic \eqref{eq:exp_cusum_recursion} is closely related to the standard CuSum \eqref{eq:cusum_recursion}, with the difference being a supplementary constant $\log \beta$ added at each step. Recently, drawing in part on the results of \cite{POOR_EXPONENTIAL}, Tartakovsky et al. \cite{tartakovsky2021optimal} used a very similar statistic for detection of transient changes, i.e. settings where the change is observable only for a limited period. An optimal procedure that maximizes the probability of detection within the transient period was derived.

\section{Recent advances in multi-stream change detection}

In modern practical applications of sequential change detection, the observed data is usually high-dimensional, observed by a large number of sensors in distinct locations. For example, when a change in the environment or radio spectrum occurs, it may affect one or multiple sensors. The problem is significantly more demanding than single-stream QCD because both the subset of affected sensors and the post-change distribution in those sensors may not be fully known or are subject to uncertainty. Additionally, the change may occur at different times in different sensors depending on the physics of the monitored phenomenon or the topology of a networked system.

\subsection{High-dimensional unknown parameters}
 
A relatively general formulation of many multi-stream change detection problems can be given as follows. The observations $\bX_n = [X_n^{(1)},...,X_n^{(K)}]$ are $K$-dimensional random vectors independent over time. However, at the same time instance the observations acquired by different sensors may be dependent. If the observations are acquired by sensors in distinct locations, it is reasonable to assume that sensor noises are independent, but the signals may be correlated. Let $f(\cdot | \btheta_n)$ denote a family of joint probability densities for $\bX_n$ parametrized by a vector $\btheta_n \in \bTheta$ of length $K$. Then $\bX_n \sim f(\cdot | \btheta_n)$, where 
\begin{equation}\label{eq:general_theta_model}
    \btheta_n = \begin{cases}
    \mathbf{0}, & n < \nu, \\
    \btheta_1, & n \geq \nu,
    \end{cases}
\end{equation}
meaning that the change manifests by altering a parameter from a known value $\bm 0$ (assumed w.l.o.g.) to an unknown $\btheta_1$. 
\subsection*{Binary $\boldsymbol\theta$}\label{sec:binary_ttheta}
We begin with the simpler setting in which the post-change model is determined solely by the subset of affected sensors. In this case, the unknown parameter $\btheta$ is naturally represented as a binary vector, where each component indicates whether the corresponding stream is affected by the change. That is, $\btheta_1 \in \bTheta =  \{0,1\}^K$, and $\btheta^{(k)} = 1$ indicates that the change affects sensor $k$. For example, if the components of $X_n$ are independent also within the same time instance and have pre- and post-change distributions $f_0^{(k)}$ and $f_1^{(k)}$, respectively, the joint density can be written
\begin{equation}\label{eq:binary_fcts}
    f(\bX_n|\btheta_n) = \prod_{k=1}^K [\theta_n^{(k)}f_1^{(k)}(X_n^{(k)}) + (1-\theta_n^{(k)})f_0^{(k)}(X_n^{(k)})].
\end{equation}

If $\btheta_1$ is known, the optimal procedure under Lorden's formulation is still the CuSum procedure with the test statistic 
\begin{equation}\label{eq:cusum_theta}
    W_n^{\btheta_1} = \max_{k\leq n}\sum_{j = k}^n \log\frac{f(X_j|\btheta_1)}{f(X_j|\mathbf{0})}
\end{equation}
When $\btheta$ is unknown, two of the most common ways of handling distributional uncertainty are to either maximize or average over $\btheta$ in \eqref{eq:cusum_theta}. The former leads to a generalized likelihood ratio (GLR) testing approach in which the unknown nuisance parameters are estimated, while the latter yields a mixture test. When the post-change parameter space is finite, the two approaches have a lot in common, as demonstrated in \cite{FELLOURIS_2016}. Indeed, the GLR-CuSum statistic is given by
\begin{align}
    &\quad \max_{k \leq n} \max_{\btheta \in \{0,1\}^K} \sum_{j=k}^n \log \frac{f(\bX_j | \btheta)}{f(\bX_j | \mathbf{0})} \\
    &=\max_{\btheta \in \{0,1\}^K} W^\btheta_n,
\end{align}
where the second line follows from interchanging the maximizations, and $W_n^\btheta$ is the CuSum statistic in which 
the set of affected sensors is $\btheta$. This procedure is equivalent to running a CuSum test for each $\btheta \in \bTheta$, and declaring a change when one of them reaches a threshold $b$. More generally, one can specify a different threshold $b_\btheta$ for each $W_n^\btheta$, $\btheta \in \bTheta$, and consider the stopping time
\begin{equation}
    T^G = \inf\left\{n: \max_{\btheta \in \{0,1\}^K} (W_n^\btheta - b_\btheta) \geq 0\right\}.
\end{equation}
If $b_\btheta$ is further chosen to be of the form $b_\btheta = b- \log p_\btheta$, where $p_\btheta > 0$ and $\sum_{\btheta \in \bTheta} p_\btheta = 1$, the stopping time can be written as
\begin{align}\label{eq:glr-binary}
        T^G &= \inf\left\{n: \max_{\btheta \in \{0,1\}^K} (W_n^\btheta + \log p_\btheta) \geq b\right\} \\
        &= \inf\left\{n: \max_{\btheta \in \{0,1\}^K} p_\btheta e^{W_n^\btheta} \geq e^b\right\}\label{eq:GLR_finite_final}.
\end{align}
Changing the maximization to a sum over $\btheta$ in \eqref{eq:GLR_finite_final}, one immediately obtains a mixture-CuSum test
\begin{equation}\label{eq:mixture-binary}
    T^M = \inf\left\{n: \sum_{\btheta \in \{0,1\}^K} p_\btheta e^{W_n^\btheta} \geq e^b\right\},
\end{equation}
where $p_\btheta$ can be interpreted as a prior distribution for $\btheta$. The GLR and mixture tests were analyzed in \cite{FELLOURIS_2016} and shown to be asymptotically optimal in a strong \emph{second-order} sense, i.e. as $\gamma \to \infty$
\begin{equation}
    \WADD(T^G) - \inf_{T : \ARL(T) \geq \gamma}\WADD(T) = O(1),
\end{equation}
for any set of affected sensors $\btheta \in \bTheta$, and similarly for $T^M$. 

The main drawback of the above procedures is that their implementation in principle requires running $|\bTheta|$ CuSum-recursions in parallel, and $|\bTheta|$ can be as large as $2^K$ if the set of affected sensors can be completely arbitrary. Significant reduction in computational complexity can be achieved by adaptive windowing \cite{FELLOURIS_2016}. A computationally simpler algorithm is the SUM-CuSum \cite{MEI_2010} for independent streams, given by
\begin{equation}
    T^\text{SUM} = \inf\left\{n : \sum_{k=1}^K W^{(k)}_n \geq b\right\},
\end{equation}
where $W^{(k)}$ is the CuSum-statistic using observations only from stream $k$ (or equivalently \eqref{eq:cusum_theta} with $\btheta = e_k$, the $k$-th basis vector). The SUM-CuSum test declares a change when the sum of local sensor CuSums crosses a threshold, making it possible to do most of the computations locally. The procedure was shown to be first-order asymptotically optimum in \cite{MEI_2010}. It also works well when the change does not necessarily occur at the same time in all affected streams. However,  unlike \eqref{eq:glr-binary} and \eqref{eq:mixture-binary}, it is not second-order asymptotically optimal.

\subsection*{Continuous $\btheta$}

In many applications, just knowing the indices of the affected sensors may not be enough to determine the full post-change distribution. For instance, the change may affect the sensors with different intensities based on e.g. the distance between the sensor and the (unknown) signal source location. For example, the sensors may be monitoring the energy of radio frequency signals in a spectrum band that attenuate inversely proportionally to squared distance. In such cases, we can still adopt the general model from \eqref{eq:general_theta_model}, but represent $\btheta \in \bTheta \subset \bbR^K$ as a general $K$-dimensional real vector. 

In principle, the general GLR and mixture-CuSum ideas from above can of course be extended to a continuous $\btheta$ without issue. However, the possibly high dimension $K$ of the unknown parameter vector $\btheta$ presents additional challenges in the continuous case. Indeed, it follows from \cite{SIEGMUND_2008} that when $f(\cdot|\btheta)$ is in the exponential family, and $\bTheta$ satisfies certain regularity conditions, the asymptotic detection delay of any procedure $T$ as ARL $\gamma \to \infty$ satisfies the following minimax lower bound:
\begin{equation}\label{eq: GLR_exponential_delay_JS}
    \sup_{\btheta \in \bTheta} I(\btheta) \WADD(T) \geq \log \gamma + \frac{K}2 \log \log \gamma + O(1),
\end{equation}
It is important to point out that in \eqref{eq: GLR_exponential_delay_JS} the second-order delay term grows linearly with $K$. For moderate $\gamma$ and large $K$, the influence of the second-order term can be greater than the first-order term.

Deriving procedures which are effective for high-dimensional $\btheta$ has attracted the most attention in the literature in the special case when $\btheta$ corresponds to the mean vector of a multivariate Gaussian distribution with diagonal covariance. 
The proposed ideas and the developed methods extend to more complicated models and different probability models. 

Xie and Siegmund \cite{XIE_2013} developed a mixture procedure where each stream is hypothesized to be affected by a change with probability $p$. The stopping time is then given by the threshold rule
\begin{equation}\label{eq:TXS}
    T^{XS} = \inf\left\{n : \max_{0\leq m \leq n}\sum_{k=1}^K \log\left(1-p + p\exp\left[\max_\btheta \sum_{j=m}^n \log\frac{f(X_j^{(k)}|\btheta)}{f(X_j^{(k)}|\bm 0)}\right]\right)\geq b\right\}.
\end{equation}
The term inside the square brackets is the GLR statistic for stream $k$ for observations between time instances $m$ and $n$. Therefore, setting $p=1$ recovers the standard GLR-CuSum. The test is shown to be effective if the post-change vector $\btheta$ is sparse, i.e. only a small proportion of its components are non-zero. In the same vein, Chan \cite{Chan2017OPTIMALDATA} proposed detectability-score transformations of stream-wise GLR statistics. The central idea is to down-weight weak local evidence that would correspond to signals too weak to be reliably detected with the current sample size. This reduces the accumulation of noise from unaffected streams and yields asymptotically near-optimal detection delay across a range of sparsity regimes. Asymptotic lower bounds for the detection delay in the limit $K\to\infty$ are also derived, which is different from the usual $\gamma \to \infty$ asymptotics in QCD. Chen et al. \cite{Chen2022High-dimensionalDetection} propose an algorithm that simultaneously tracks various likelihood ratios corresponding to individual or multiple sensors, and uses hard thresholding to suppress likelihood ratios likely corresponding to noise in the sparse case. Using different levels for hard threshold in parallel makes the test adaptable to different levels of sparsity in the post-change mean vector $\btheta$.

The usefulness of thresholding in detecting a change to a sparse $\btheta$ is related to the effectiveness of thresholding in sparse estimation, see e.g. \cite{candes2006modern} and references therein. Indeed, a change detection problem with unknown parameters may be seen as a combination of an estimation and detection problem; the better $\btheta$ is known (estimated) the easier a change is to detect. Another popular technique used in high-dimensional problems is shrinkage estimation \cite{FOURDRINIER_BOOK}. A natural way to incorporate a particular estimator to a change detection algorithm is to use a so-called adaptive change detection procedure \cite{LORDEN_2005, LORDEN_2008_COMPUTATIONALLY, CAO_2018, XIE_2023}. For example, in the CuSum recursion \eqref{eq:cusum_recursion}, one can replace the unknown $f_1$ (here $f(\cdot|\btheta)$) with a plug-in $f(\bX_n|\bthetahat_n)$, where $\bthetahat_n$ is some estimator of $\btheta$ based on the data up to time instance $n$. Xie et al. \cite{XIE_2023} proposed estimating $\btheta$ based on a sliding window of past $w$ observations, i.e. $\bthetahat_n = \bthetahat(\bX_{n-1},\dots, \bX_{n-w})$. This yields the window-limited CuSum statistic
\begin{equation}\label{eq:WL_stat}
    W^{\text{WL}}_n = \max\{0, W_{n-1}\} + \log\frac{f(\bX_n|\bthetahat_n)}{f({\bX_n | \bm 0)}}, \quad W_0 =0.
\end{equation}

An important feature of this statistic is that it can be calculated recursively, which is not in general the case for e.g. GLR type tests, as remarked in Section~\ref{sec:binary_ttheta}. Another advantage is that prior knowledge/structure (such as sparsity) about $\btheta$ can be easily incorporated by selecting the estimator $\bthetahat$ accordingly. For example, shrinkage estimators reduce estimator variance by introducing some bias, often leading to reduced loss under an appropriate loss function. Wang and Mei \cite{WANG_2015} considered thresholding as well as linear shrinkage estimators of the form $a\thetaML + b$, where $\thetaML$ is the maximum likelihood estimator. Later, it was demonstrated in \cite{Halme2025QuickestEstimator} that using the James-Stein (JS) estimator \cite{JAMES_STEIN_1961} in place of the standard maximum likelihood estimator in the WL-CuSum test \eqref{eq:WL_stat} can reduce detection delay for all values of $\btheta$ in the Gaussian case. This stems from the fact that the James-Stein estimator dominates the ML estimator in terms of mean-square error when $K \geq 3$ \cite{JAMES_STEIN_1961}. The JS-estimator is most effective when $\btheta$ is dense rather than sparse, i.e. when the change event causes subtle disturbances to many sensors. A simple way to create a procedure that is effective for both dense and sparse $\btheta$ simultaneously would be to run two versions of \eqref{eq:WL_stat} in parallel: one with an estimator $\bthetahat_n$ effective for dense changes (such as JS) and one for sparse changes (such as a hard-thresholding estimator \cite[Ch. 5]{candes2006modern}). Forming and analyzing such combined procedures in a principled manner is an interesting topic for future research. 

\subsection*{Numerical comparison of methods for high-dimensional Gaussian mean-shifts}

As an illustration of the properties of the methods discussed above, we present a brief numerical comparison. We consider a model with $K=100$ data streams, where before the change
\[
f_0=\calN(\bm 0,I_K),
\]
and after the change
\[
f_1=\calN(\btheta,I_K).
\]
Although simple, this setting is not unduly restrictive: more general Gaussian models with correlated components can often be reduced to this form by whitening, provided that the pre-change covariance matrix is known or can be estimated accurately.

We compare six procedures: $T^{XS}$ in \eqref{eq:TXS} from \cite{XIE_2013}, Chan's detectability score procedure \cite{Chan2017OPTIMALDATA}, the OCD procedure from \cite{Chen2022High-dimensionalDetection}, two versions of the WL-CuSum procedure \eqref{eq:WL_stat} from \cite{XIE_2023} based on the ML and James--Stein estimators \cite{Halme2025QuickestEstimator}, and the classical GLR-CuSum procedure \cite{LAI}. Since XS, Chan, and GLR do not admit recursive updates, we implement window-limited versions with window size 200. For the XS and Chan procedures, we set $p=1/\sqrt{K}=0.1$ as recommended in \cite{Chan2017OPTIMALDATA}. Parallel versions of the WL-CuSum tests are run (see \cite{XIE_2023}) with window-sizes ranging from 30 to 80. The XS, Chan and OCD procedures are implemented using the R-package \texttt{ocd} \cite{CHEN_DATA}. The stopping thresholds of all procedures are set via 1000 Monte-Carlo simulations such that $\ARL(T) \approx 1000$.

\begin{figure}
    \centering
    \includegraphics[width=0.9\linewidth]{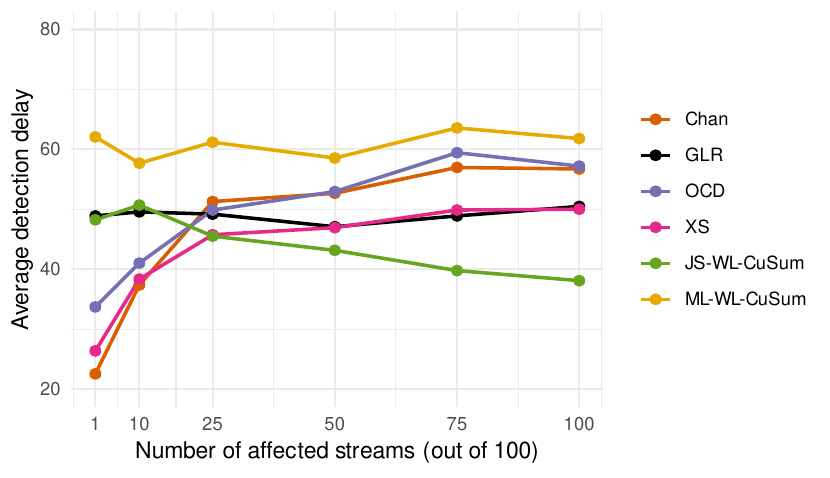}
    \caption{Average detection delays highlight a clear dependence on the sparsity of the change. For very sparse changes, XS \cite{XIE_2013}, Chan \cite{Chan2017OPTIMALDATA}, and OCD \cite{Chen2022High-dimensionalDetection} procedures achieve the smallest delays, whereas their advantage diminishes as more streams are affected. In contrast, the James-Stein WL-CuSum \cite{Halme2025QuickestEstimator} iis ncreasingly effective in denser regimes, while the ML-WL-CuSum and GLR procedures remain comparatively stable across sparsity levels. Overall, it is observed that substantial gains over standard methods are possible in high-dimensional problems when the detection procedure is matched to the structure of the post-change signal.}
    \label{fig:placeholder}
\end{figure}

We consider changes affecting $L\in\{1,25,50,100\}$ of the $K=100$ streams. For a given value of $L$, the post-change vector $\btheta$ is generated at every MC iteration as follows: for $1\leq k\leq K$ we first sample $\xi^{(k)} \sim \calN(0.5, 1)\ind{k\leq L}$ and then set $\btheta^{(k)} = \xi^{(k)}/\sqrt{\sum(\xi^{(k)})^2}$. This normalization ensures that the KL-divergence satisfies $I = \norm \btheta^2/2=1/2$ at every iteration. The change-point is fixed at $\nu=100$. Figure~\ref{fig:placeholder} reports Monte Carlo estimates of the average detection delay $\Ex_\nu(T-\nu\mid T\ge \nu)$. The results show a clear dependence on sparsity. For sparse changes, the XS, Chan, and OCD procedures perform best. As the change becomes denser, their advantage diminishes, whereas the James--Stein WL-CuSum becomes increasingly competitive. By contrast, the ML-WL-CuSum and GLR procedures are comparatively insensitive to the number of affected streams. Importantly, for both very sparse and very dense $\btheta$, substantial performance gains over standard methods like GLR are available when the number of sensors is large.

\subsection{QCD under sampling constraints}

In the previous subsection, we assumed that all sensors acquire data continuously at each time instance. In practice, it is often necessary to reduce the communication and sampling overhead by sampling only some subset of streams at a given time. By designing intelligent sensing policies and methods for sensor selection, significant savings in computational complexity and communications overhead may be achieved.  Let us denote by $\mathcal{A}_n \subset \{1,...,K\}$ the indices of streams monitored at time $n$, meaning that with some abuse of notation
\begin{equation}
    X_n^{(k)} = \begin{cases}
        X_n^{(k)} \quad &\text{if }k \in \mathcal{A}_n \\
        \emptyset &\text{otherwise},
    \end{cases}
\end{equation}
where $\emptyset$ denotes an empty symbol/missing observation. The sampling choice $\mathcal{A}_n$ can depend on previous data, i.e. $\mathcal{A}_n \in \calF_{n-1}$. The size $\mathcal{A}_n$ is usually constrained, for example $|\mathcal{A}_n| \leq L~\forall n$ if at most $L < K$ of the sensors can be deployed at a given time. Under such a constraint, the typical QCD objective is now expanded to designing both a sampling policy and a detection rule. 

Dragalin \cite{Dragalin1996ASystem} considers the case $L = 1$ and proposes a scanning rule for stream selection in which a particular stream is sampled until its CuSum statistic hits either $0$ or a stopping threshold $b$. If the statistic hits $0$, the sampling moves to the next data stream. The procedure is shown to be asymptotically optimum under Pollak's CADD. Lai et al. \cite{Lai_quickest_search} consider a closely related problem with $L = 1$ in a sequential hypothesis testing context. In an infinite sequence of streams, a signal is present independently in each stream with equal probability. The objective is to find an affected stream as quickly as possible while sampling only a single stream at a time. It is shown that within a Bayesian formulation seeking minimal delay and false detection probability, successively running a CuSum test until either reset and switch to another stream, or detection in the current stream is exactly optimal. Recently, \cite{Xu_sampling_control} considered a similar rule in a QCD context, and showed that such a cyclic procedure is asymptotically second-order optimal under Lorden's WADD. Moreover, the data distributions in different streams are allowed to be different. 

The general idea in the above works may be extended to the case when $1 \leq L \leq K$ streams are observed at a time by e.g. computing and updating a detection statistic for each stream, and sampling those with the $L$ currently highest values. However, as observed in \cite{LIU_MEI}, especially if $K$ is much larger than $L$, it may be beneficial to ``increase exploration'' by incrementing the detection statistics of all unobserved streams at each step by a small constant $\Delta$. This way, the procedure cannot get stuck sampling a single stream for a long time. Determining an ``optimal'' value for $\Delta$ is not straightforward, and even the analysis of the procedure is challenging for $\Delta > 0$, but small values of $\Delta > 0$ are shown to perform better empirically than $\Delta = 0$ \cite{Xu2025ImpactDetection}. 

The above works assume that sampled data streams are independent, or at the very least that the change is observable from the marginal distributions of individual streams. If the data are acquired by distinct sensors, noise can be assumed to be independent while the signals may be correlated or otherwise dependent. More general observation models were studied recently in \cite{Chaudhuri2024RoundData}, where the change is assumed to affect the joint distribution of an unknown subset of sensors; for example the correlation structure. Conditions under which the cyclic, rotating order sampling policy remains asymptotically optimal were derived. Previously in e.g. \cite{Xu_sampling_control}, asymptotic optimality when only a single sensor observes the change was established with respect to a procedure that has access to all samples. In \cite{Chaudhuri2024RoundData} asymptotic optimality is studied only with respect to procedures that satisfy the sampling constraint. In this class, the simple sampling policy of \cite{Xu_sampling_control} is first-order asymptotically optimal even when an arbitrary number of sensors observe the change under the assumption of homogeneity of the statistical models.

The cyclic policies are computationally simple to implement and don't require storing a lot of past data. From a broader perspective, designing the sampling policy can be seen to be closely related to the multi-armed bandit (MAB) \cite{Lattimore_Szepesvari_2020} problem. The core challenge in the MAB problem is the exploration vs. exploitation tradeoff: choosing the arm (stream) with the highest expected reward (exploitation) vs. choosing an arm (stream) that may have a higher reward than the current best arm (exploration). As an example, Restless Multiarm Bandit (RMAB) \cite{Oksanen2015AnNetworks, Lunden2015SpectrumAdvances}
and Bayesian Thompson sampling methods have been used for learning sensing policy and stream selection while balancing between exploitation and exploration in multisensor systems. In a change detection setting, however, the objective is not to minimize e.g. cumulative regret, but for instance average detection delay.  In case of long-term rewards, associating a per-action reward to each action (sampling choice) is therefore not obvious, although a large observed likelihood ratio could represent a reasonable proxy. In fact, from this perspective, the cyclic sampling policies discussed above can be seen as greedy (or follow-the-leader) selection policies. In general, myopic learning policies that maximize immediate or near-term rewards are more appropriate for stream or sensor selection in the QCD context. 

A natural and fundamental algorithm for sampling that often works well in practice is Thompson sampling in a Bayesian framework. The idea is to define a distribution over possible environments, and at each time make the optimal action with respect to an environment drawn at random from the current posterior distribution. In a QCD context, this would mean specifying a (joint) prior distribution for the change-points, computing the posterior probabilities of change having happened at each stream (i.e. the Shiryaev statistics \eqref{eq:shiryaev_time}), and sampling according to these probabilities. Such a policy is studied in \cite{Nitzan2020BayesianCommunication}. 
If specifying a prior distribution is infeasible, using the Shiryaev-Roberts statistics \eqref{eq:SR_stat} is a logical substitute. A Thompson sampling related sampling policy based on randomized SR statistics is proposed and analyzed in \cite{Zhang2023BanditControl}.

Finally, Veeravalli et al. \cite{VVV_Controlled} consider a generalized setup where possible sampling actions $\calA_n$ do not have to correspond to subsets of sensors, but represent general (finite) set of sampling actions. For example, different sampling actions could represent different ways of steering a sensor array. 
Depending on the true target parameters $\btheta$, some sampling actions $a$ are more effective than others for distinguishing the change. The effectiveness is measured by the post-change KL-divergence 
\begin{equation}
    I^a(\btheta) = \int f^a(x|\btheta) \log \frac{f^a(x|\btheta)}{f^a(x|\mathbf{0})},
\end{equation}
where $f^a(\cdot|\btheta)$ is the sampling distribution of the data under action $a$ when the true parameter is $\btheta$. For example, if $\btheta$ is a binary vector as in \eqref{eq:binary_fcts}, $\mathcal{A}_n$ is a subset of sensors, and the components of $\bX_n$ are independent, 
\begin{equation}
    f^{\mathcal{A}_n}(\bX_n|\btheta) =  \prod_{k \in \mathcal{A}_n} [\theta_n^{(k)}f_1^{(k)}(X_n^{(k)}) + (1-\theta_n^{(k)})f_0^{(k)}(X_n^{(k)})].
\end{equation}
Naturally, the optimal action is $a^* = {\arg \max}_a I^a(\btheta)$, but this depends on the unknown $\btheta$. As a solution, it is proposed in \cite{VVV_Controlled} to 1) estimate $\btheta$ with a maximum likelihood sliding window estimator $\bthetahat_n$ based on past observations, and 2) select the next action $\mathcal{A}_n$ as either $\mathcal{A}_n = {\arg \max}_a I^a(\bthetahat_n)$ or uniformly at random at deterministic exploration times. The form of the sampling policy goes all the way back to Chernoff and the sequential design of experiments in hypothesis testing \cite{Chernoff_sequential}. Sampling policies based on Chernoff's rule are first-order asymptotically optimum across a wide range of controlled sensing tasks. However, its (empirical and theoretical) performance can be suboptimal when the choice of an action can change the effect of future actions, as noted recently in \cite{Tajer2022ActiveNetworks}. 

\subsection{QCD with False Discovery Rate control}

So far, we have considered settings where only a single change affecting the monitored system is detected. In many applications changes can happen at different times at different sensors, possibly caused by the same or different phenomena. For example, in a sensor network, each sensor may be monitoring a different physical process, and a change in one process may or may not affect the others. When a change takes place in one stream, it can require an action to be taken only in that stream. In these cases it may be necessary to detect a change-point separately in each stream. Moreover, it is of interest to consider performance criteria that measure the performance across the streams as a whole, instead of separately optimizing the performance in each individual stream. For example, in a cognitive radio network \cite{LAI_GLOBECOM}, secondary users opportunistically sense the channels for available transmission opportunities. Even if the individual channels are independent, it is important to limit the total amount of interference or collisions caused by the SUs as a whole. 

Let us denote the time of detection in the $k$th stream by $T^{(k)}$, and by $\nu^{(k)}$ the actual change-point in stream $k$. If a QCD procedure that controls the ARL \eqref{eq: arl_definition} individually is applied to each stream, one approximately obtains control over the average \emph{number} of false-alarm events that can occur over all streams in a given time interval. In practice, it might be more useful to have control over the \emph{proportion} of false detections among all detections. To this end, a natural criterion well-suited for large-scale applications is the false discovery rate (FDR) \cite{BENJAMINI_HOCHBERG}. In a change-detection setting, the FDR at time $n$ for a set of stopping times $\calT = \{T^{(1)},\dots,T^{(K)}\}$ is defined as the expected ratio of false detections $V(n)$ made by time $n$ to all detections $R(n)$ by $n$:
\begin{equation}\label{eq:fdr_def}
\FDR_\bnu(\calT, n) = \Ex_\bnu\left[\frac{V(n)}{\max\{1,R(n)\}}\right]
\end{equation}
where
\begin{align}
    V(n) & = \sum_{k=1}^K\ind{T^{(k)} < \nuk, T^{(k)}\leq n}    \\
    R(n) & = \sum_{k=1}^K\ind{T^{(k)} \leq n}.
\end{align}
The denominator in \eqref{eq:fdr_def} includes the maximum only to prevent division by zero.

In the FDR literature, the typical starting point is a set of p-values corresponding to fixed-sample size hypothesis tests and the ubiquitous Benjamini-Hochberg procedure \cite{BENJAMINI_HOCHBERG}. The BH-procedure was extended to sequential hypothesis tests by Bartroff and Song \cite{Bartroff2020SequentialRates}. Utilizing some ideas from \cite{Bartroff2020SequentialRates}, Chen et al. \cite{CHEN_ZHANG_POOR_BAYESIAN_JOURNAL} were the first to study multi-stream QCD procedures with FDR control. Under the Bayesian change-point formulation (Sec.~\ref{sec:bayes_object}), a procedure that utilizes the posterior probabilities \eqref{eq:posterior_prob} as the test statistics in a sequential BH-type procedure is proposed. It is shown that the procedure controls the FDR averaged over the distribution of the change-points $\pi(\bnu)$ asymptotically as $K \to \infty$, i.e. for all fixed $n$ and a user-defined FDR-threshold $0 < \alpha < 1$
\begin{equation}\label{eq:asym_fdr}
    \FDR^\pi(\calT, n)= \sum_\bnu \pi(\bnu) \FDR_\bnu(\calT, n) \leq \alpha + o(1).
\end{equation}
Later, it was observed in \cite{Nitzan2020BayesianCommunication} that in the Bayesian setting it is in fact sufficient to use a single common stopping threshold for all streamwise tests, and such a procedure controls the FDR non-asymptotically. That is, \eqref{eq:asym_fdr} is achieved without the $o(1)$ term. In a simple setting, even an exactly optimal procedure for $\FDR^\pi(\calT, \infty)$ can be derived \cite{asilomar_2023} as a concatenation of Shiryaev tests \eqref{eq:shiryaev_time} with specifically chosen stopping thresholds.

Under a non-Bayesian formulation, the natural objective in the spirit of Pollak's CADD \eqref{eq:CADD_def} would be to control the worst-case FDR over all values of the change-point, namely
\begin{equation}\label{eq:sup_fdr}
    \sup_\bnu \FDR_\bnu(\calT,n) \leq \alpha
\end{equation}
This condition is clearly stronger than in the Bayesian setting, since controlling the FDR for all individual $\bnu$ clearly implies that the FDR averaged over any distribution $\bnu \sim \pi$ is controlled as well. For example, \eqref{eq:sup_fdr} implies that $\FDR_\infty(\calT, \infty) \leq \alpha$, which in turn implies
\begin{equation}
    P_\infty\left(\bigcup_{k = 1}^K \{T^{(k)} < \infty \}\right) \leq \alpha.
\end{equation}
That is, in the case where no change ever occurs, the probability of any procedure declaring a change should be smaller than $\alpha$. Such a condition is more stringent than the typical ARL used in single-stream change detection. Indeed, $P_\infty(T< \infty) < 1$ implies $\ARL(T) = \infty$. It is expected that a procedure that controls \eqref{eq:sup_fdr} will suffer large detection delays, at least in the worst-case. Nonetheless, procedures that control \eqref{eq:sup_fdr} were developed in \cite{CHEN_ZHANG_POOR_NON_BAYESIAN, CISS2020}. The methods consist of running a CuSum test in each stream with certain stopping thresholds that are increasing over time. As a result, the average detection delay can be large if the change happens late, since by the time the change happens the stopping threshold has increased to a high level. 

Recently, \cite{DANDAPANTHULA} consider a slightly more general version of \eqref{eq:sup_fdr}, where FDR is controlled not only at any fixed $n$, but also any stopping time $\tau$. For example, $\tau$ could be the time $m$ changes have been detected. Since constructing a detection procedure with small (or even finite) worst-case average detection delay is infeasible due to the reasons discussed above, a criterion called error over patience (EOP) is proposed:
\begin{equation}
    \text{EOP}(\calT) = \sup_{\bnu \geq 1}\sup_{\tau} \frac{\Ex(\FDR_\bnu(\calT,\tau))}{\Ex(\tau)},
\end{equation}
where the inner supremum is taken over all stopping times (Definition~\ref{def:stoppingtime}).
In essence, the criterion allows the FDR to be larger if the process is terminated late, i.e. when $\tau$ is large. This way, control over the FDR does not need to be maintained for an indefinitely long time.  

The definition of FDR we have adopted \eqref{eq:fdr_def} measures the quality of detections made up to time $n$. In some applications, it could be of interest to consider the ``quality'' of the streams still active at time $n$. For example, the false non-discovery rate (FNR) at time $n$ defined as 
\begin{equation}\label{eq:fnr_def}
\text{FNR}_\bnu(\calT, n) = \Ex_\bnu\left[\dfrac{\sum_{k=1}^K {\ind{\nu^{(k)} \leq n, T^{(k)} > n}}}{\max\{1,K-R(n)\}}\right]
\end{equation}
measures the proportion of post-change streams among all streams where a change has not yet been detected. Discussions of many related performance metrics such as FNR and algorithms controlling them can be found in \cite{Chen2023CompoundStreams} and \cite{Lu2022OptimalMeasures}.

\section{Future outlook}

In this brief overview, the discussion was mostly limited to parametric models. In modern applications, it may not always be possible to formulate the data distributions using parametric models. Many approaches that replace the likelihood ratio with some non-parametric data-driven statistic and empirical probability models have been proposed recently. For example, tests based on kernel minimum mean discrepancy (MMD) have been recently proposed \cite{LI_SCAN, Wei2026OnlineDetection} as interesting ideas.  Many other ways of finding statistics that discriminate between the pre-change and post-change distributions in place of the likelihood ratio naturally exist as well, see for example \cite[Sec. IV-B]{XIE_2021_REVIEW}. Empirical likelihoods (EL) and likelihood ratios provide a data driven, non-parametric approach while being able to take advantage of some benefits of likelihood-based sequential inference. Moreover, conformal p-values computed using a training data from the pre-change probability model to indicate a change in distribution, outliers or out of distribution data may be employed. They have found applications in sequential inference in single and multi-stream settings and multisensor systems performing spatio-temporal inference. 

One recently proposed general framework for how statistics other than the likelihood ratio can be used for (single-stream) QCD is proposed by Shin et al. \cite{Shin2023E-detectors:Detection}. It is shown that the general CuSum and Shiryaev-Roberts procedures can also be used when the likelihood ratio is replaced by a so-called e-process. E-processes may be seen as a nonparametric generalization of the likelihood ratio. 

Optimal detectors and estimators are derived based on a particular parametric signal model and strict assumptions on the underlying probability models. However, optimality is only achieved when the underlying assumptions hold. The performance of optimal procedures may deteriorate significantly, even for minor departures from the assumed model. 
Robust statistical methods are semi-parametric, namely take advantage of parametric models but trade off optimality to reliability.  Statistically robust detectors are close-to optimal in nominal conditions and highly reliable for real-life data, even if the assumptions are only approximately valid. Robust procedures are an attractive approach in QCD problems where outliers may be present in the data, or the assumptions in probability models do not hold exactly, see e.g. \cite{UNNI_ROBUST, XIE_ROBUST, Hou2025RobustControl} and references therein. 
Robust procedures can distinguish among data from different populations and outliers. The parameters of the distribution may then be estimated from the population associated with the majority of the data, and the outliers and data from minority population will have a bounded influence. 

Machine learning  approaches provide an attractive approach for change-point detection in the presence of modeling deficits. Probability models and their parameters can then be learned from data. In highly dynamic and high-dimensional scenarios, learning sensing policies allows for allocating sensing, computation, and inference resources to a subset of data streams where the change is likely to take place soon. Consequently, significant savings in computation and energy consumption may be achieved. Different variants and recent advances of reinforcement learning, including MAB and Partially Observable Markov Decision Processes (POMDP) provide a framework for learning optimized sensing policies. Statistical censoring methods can be used to reduce data transmissions and energy consumption in networked, multi-stream change-point detection since only informative decision statistics are transmitted. Only minor performance losses in inference are experienced even under rather strict communications rate constraints. 

No matter the exact approach that is used in a non-parametric setting, the general challenges highlighted in this paper still persist. High dimensionality, sensing constraints, and multiplicity of decisions continue to complicate both algorithm design and performance analysis.

In a broader context, detecting rapid changes, disruptions and adversarial events is a key capability in resilient systems that are able to identify such undesirable events, recover and reconfigure to continue normal operation, and provide desired trustworthy sequential inference performance with rigorous performance guarantees. Detecting changes is also crucial in creating situational awareness exploited by any adaptive or cognitive system.

\section{Conclusion}

This overview surveyed some key models and methods for multi-stream quickest change detection, with emphasis on three recurring challenges: high-dimensional unknown post-change parameters, multiple-hypothesis testing inspired error rate control, and adaptive sensing under resource constraints. Across these settings, a common message is that effective procedures must balance statistical optimality with computational and operational feasibility.

Several important challenges remain open. In high-dimensional settings, a central difficulty is to design procedures that retain strong optimality properties while remaining computationally feasible when the number of streams is very large. Beyond the parametric models emphasized here, there is also a growing need for methods that operate reliably when the pre- and post-change distributions are only partially specified or must be learned from data. In controlled sensing problems, an important direction is to develop adaptive sampling rules with rigorous performance guarantees under dependence, heterogeneity, and resource constraints. Finally, in settings with multiple streamwise detections, further work is needed to combine local detection efficiency with performance criteria that measure the global performance. A related challenge is to formulate performance criteria that are both achievable in practice and sufficiently stringent to control the error events and tradeoffs that matter most in applications.

\vspace{6pt} 



\subsubsection*{Acknowledgement}
The authors would like to thank Dr. H. Vincent Poor for insightful discussions on the topic of quickest change detection, and for the fruitful long-term collaboration  on a variety of signal processing topics.

\bibliographystyle{ieeetr}
\bibliography{references}

\end{document}

%% file: shorts.tex
\def\Xseq{{X_1,X_2,\ldots}}

\def\Nmax{{N_{\max}}}

\def\PFA{{\textsc{PFA}}}
\def\ADD{{\textsc{ADD}}}
\def\ARL{{\textsc{ARL}}}
\def\FDR{\textsc{FDR}}
\def\WADD{{\textsc{WADD}}}
\def\CADD{{\textsc{CADD}}}

\def\thetahatnorm{{\lVert \bthetahat_{n-1} \rVert^2}}
\def\thetanorm{{\lVert \btheta \rVert^2}}
\def\JSplusnorm{{\lVert \thetaJSplus \rVert^2}}

\def\llr{{\ell(\bX_n, \hat\btheta_{n-1})}}
\def\llrMLE{{\ell(\bX_n, \hat\btheta^\textsc{ML}_{n-1})}}
\def\llrJS{{\ell(\bX_n, \hat\btheta^{\textsc{JS}^+}_{n-1})}}

\def\JHAT{{\mathsf J(\bthetahat)}}
\def\IHAT{{\mathsf I(\bthetahat)}}
\def\I{{\mathsf I}_1}

\def\IMLE{{\mathsf I\left(\thetaML\right)}}
\def\IJS{{\mathsf I\left(\thetaJSplus\right)}}
\def\JJS{{\mathsf J\left(\thetaJSplus\right)}}
\def\JMLE{{\mathsf J\left(\thetaML\right)}}

\def\bmu{{\boldsymbol{\mu}}}
\def\bI{{\mathbf{I}}}
\def\bZ{{\mathbf{Z}}}
\def\bv{{\mathbf{v}}}
\def\bX{{\mathbf{X}}}
\def\bnu{{\boldsymbol{\nu}}}

\def\Xk{X^{(k)}}
\def\nuk{{\nu^{(k)}}}

\def\JSplus{{\mathsf{JS}^+}}

\def\thetaML{{\bthetahat^{\normalfont{\textsc{ML}}}}}
\def\thetaJS{{\bthetahat^{\normalfont{\textsc{JS}}}}}
\def\thetaJSm{{\bthetahat^{\textsc{JS}}_m}}
\def\thetaJSmplus{{\bthetahat^{\textsc{JS}}_{m+}}}
\def\thetaJSVplus{{\bthetahat^{\textsc{JS+}}_\mathbb{V}}}
\def\thetaJSV{{\bthetahat^{\textsc{JS}}_\mathbb{V}}}

\def\thetaJSplus{{\bthetahat^{\normalfont{\textsc{JS}}}_+}}

\def\x{{\bm x}}
\def\barX{{\overbar{\bX}}}

\def\realR{{\mathbb R}}

\def\mfR{{\mathfrak{R}}}

\def\calA{{\cal{A}}}
\def\calL{{\cal L}}
\def\calN{{\cal N}}
\def\calF{{\cal F}}
\def\calC{{\cal C}}
\def\calD{{\cal D}}
\def\calI{{\cal I}}
\def\calJ{{\cal J}}
\def\calT{{\cal T}}
\def\calH{{\cal H }}
\def\Ex{{\mathbb{E}}}
\def\P{{\mathbb P}}
\def\Var{{\mathrm{Var}}}
\def\bI{{\mathbf{I}}}
\def\btheta{{\boldsymbol{\theta}}}
\def\bthetahat{{\skew{3}{\hat}{\btheta}}}

\def\bTheta{{\boldsymbol{\Theta} }}

\def\bbR{{\mathbb{R}}}
\def\bbV{{\mathbb{V}}}
\def\projV{{P_\mathbb{V}}}

\def\DKL{{D(f_1 || f_0)}}

\def\bv{{\mathbf{v}}}
\def\bX{{\mathbf{X}}}
\def\bY{{\mathbf{Y}}}
\def\bZ{{\mathbf{Z}}}

\def\ddX{{\frac{\partial}{\partial \bar{\bX}}}}

\def\WGLR{{W_n^{\scriptscriptstyle{\text{GLR}}}}}
\def\TGLR{{T^{\textsc{G}}_b}}
\def\WSRRS{{W_n^{\scriptscriptstyle{\text{SRRS}}}}}
\def\TSRRS{{T^\textnormal{{\textsc{S}}}_b}}
\def\TSRRSgamma{{T^{\textsc{S}}_\gamma}}
\def\TWL{{T^{\normalfont{\textsc{WL}}}_b}}
\def\TSR{{T^\textsc{SR}_b}}
\def\TCUSUM{{T^\textsc{C}_b}}

\newcommand{\ind}[1]{\textbf{1}{\left\{#1\right\}}}
\newcommand\given[1][]{\:#1\vert\:}
\newcommand\norm[1]{\lVert#1\rVert}
\newcommand\delay[1]{\calD^\btheta\left(#1\right)}
\newcommand\MSE[1]{\mathsf{MSE}_\btheta(#1)}
\newcommand{\KL}[2]{D(#1 \,\|\, #2)}
\def\define{\coloneq}

\newcommand{\overbar}[1]{\mkern 1.5mu\overline{\mkern-1.5mu#1\mkern-1.5mu}\mkern 1.5mu}